\documentclass[12pt, a4paper]{amsart}

\usepackage{amsmath,amssymb,amscd,amsfonts,verbatim}

\newtheorem{thm}{Theorem}[section]
\newtheorem{lem}[thm]{Lemma}

\newtheorem{prop}[thm]{Proposition}
\newtheorem{cor}[thm]{Corollary}

\newtheorem{ex}[thm]{Example}

\newcommand\si{{\sigma}}
\newcommand\cF{{\mathbb{F}}}

\newcommand{\Pic}{\operatorname{Pic}}
\newcommand{\rank}{\operatorname{rank}}
\newcommand{\Aut}{\operatorname{Aut}}

\newcommand\cod{{\text{\rm codim}}}

\newcommand\alb{\text{\rm a}}

\newcommand\Alb{\text{\rm A}}
\newcommand\lra{\longrightarrow}

\newcommand\ot{{\otimes}}
\newcommand\OO{{\mathcal{O}}}

\newcommand\pp{{\mathbb{P}}}
\newcommand\ox{{\omega _X}}

\newcommand\FF{{\mathcal{F}}}
\newcommand\F{{\mathcal{F}}}

\newcommand\Z{{\mathbb{Z}}}

\newcommand{\Si}{\Sigma}

\numberwithin{equation}{section}

\title[Products of curves of genus $2$]
{ Birational characterization of products of curves of
genus 2}
\author {Christopher D. Hacon and Rita Pardini}
\date{}
\begin{document}
\abstract Let $X$ be a variety of maximal Albanese dimension.
In this paper we prove that
if $\chi (\omega _X )=1$ then $q(X)\leq 2\dim X$ and if $q(X)=2\dim X$, 
then $X$ is birational to a product of
curves of genus $2$.\newline
\endabstract
\maketitle

\section{ Introduction}
The purpose of this paper is to show how a result of Beauville 
in the classification
theory of irregular surfaces of general type can be generalized to 
higher dimensional varieties
of maximal Albanese dimension. 
 
It is well known that for  a surface $X$
of general type one has  $\chi (\omega _X)>0$. Therefore surfaces 
with $\chi (\omega _X )=1$ can be regarded as the ``first case in 
the classification'' and
for this reason they have always been of interest to people working on the 
theory of surfaces.
The case in which the irregularity $q(X)$ is equal to $0$ has been 
investigated for a long time
  by several authors, but a complete classification seems at the moment out 
of reach, due to the
existence of a great number of different examples.
On the other hand, the following result of Beauville (\cite{Be})  
suggests that the
classification should be easier when the irregularity $q(X)$  is large:

{\it 
If $X$ is a surface
of general type with $\chi (\omega _X)=1$ then $q(X)\leq 4$
and if $q(X)=4$, then $X$ is birational to the product of
two curves of genus $2$.} 

In fact, more recently, Hacon and Pardini \cite{HP1}
and Pirola \cite{Pi} have classified surfaces
of general type with $\chi (\omega _X)=1$ and $q(X)=3$, showing that they belong to two families.
Surfaces of general type with $\chi (\omega _X)=1$ and $q(X)=1,2$
have proven to be much harder to understand and they
are still an active topic of research.

It is a natural question to try and generalize
the above results to higher dimension. In general, if $\dim X\geq 3$,
one only expects results of a qualitative nature. Indeed,   in Section
\ref{Sexamples} we give a construction of threefolds of general type  
with $\chi(\omega _X)=1$ and
arbitrarily large irregularity, thus showing that it is not possible 
to generalize Beauville's
statement verbatim. Therefore, the situation appears to be substantially 
more complicated than
that of surfaces. However, Green and Lazarsfeld have made the surprising 
remark that varieties of
maximal Albanese dimension seem to behave similarly to
surfaces. For example,
they prove the inequality  $\chi (\omega _X)\geq 0$ for a variety $X$ of 
maximal Albanese
dimension (cf.
\cite{GL1}). Another example of this analogy is the fact, proven by Chen 
and Hacon (\cite{ChH1}),
that for  a variety of general type and 
maximal Albanese dimension, the $6-$canonical map is
always birational.

In this note we give  an additional instance of the analogy between  
surfaces 
and varieties of maximal Albanese dimension,  by proving  the following
generalization of the above mentioned result of Beauville:
\medskip

\noindent{\bf Theorem.}
{\em Let $X$ be a smooth projective variety of 
maximal Albanese dimension. If $\chi (\omega _X)=1$, then $q(X)\leq
2\dim X$ and if $q(X)=2\dim X$, then $X$ is birational to a product 
of curves of genus $2$.}
\medskip

The proof of the Theorem above
relies on the generic vanishing theorems of Green and Lazarsfeld
\cite{GL1}, \cite{GL2} and their generalizations \cite{Si}, \cite{ClH},
\cite{ChH2}, \cite{Hac}.
It would be interesting to see if these techniques can be successfully applied
in the case when $\dim X<q(X)<2\dim X$. However, already the case
$\chi(\ox)=1$ and $q(X)=2\dim X-1$ seems considerably harder. As we have already mentioned,
there is a complete classification in dimension 2, but the number of examples grows with the
dimension and one does not have  clear picture of the situation.
\medskip

\noindent{\bf Acknowledgments.} The first author was partially supported
by NSA research grant no: MDA904-03-1-0101 and by a grant from the
Sloan Foundation. The second author was partially supported by  P.R.I.N. 2002
``Geometria sulle variet\`a algebriche'' of Italian M.I.U.R. and is a a member of
G.N.S.A.G.A. of C.N.R.

\medskip
\noindent{\bf Notation and conventions.}
We work over the field of complex numbers; all varieties are projective.
We identify Cartier divisors and line bundles on
a smooth variety, and we use the additive and
multiplicative notation interchangeably. If
$X$ is a smooth projective variety, we let $K_X$ be a  canonical divisor,
so that
$\omega _X=\OO_X(K_X)$, we denote by
$\kappa(X)$ the Kodaira dimension and  by
$q(X):=h^1(\OO_X)=h^0(\Omega^1_X)$ the {\em irregularity}. 
We denote by  $\alb \colon X\to \Alb(X)$ the Albanese map and by
$\Pic^0(X)$ the dual abelian variety to $\Alb (X)$, which parameterizes all
topologically trivial line bundles on $X$. $\Pic^0(X)_{tors}$
 denotes  the set of torsion points in $\Pic^0(X)$.
If the Albanese morphism is generically finite onto its image,
we say that $X$ is of {\it maximal Albanese dimension}.
If $\kappa (X)=\dim X$, we say that $X$ is of {\it general type}.
\section{Preliminaries.}\label{pre}

The main ingredients in the proof  of Theorem \ref{mainT} are the relative vanishing
theorems due  Koll\'ar and a generalized version of the generic vanishing theorems of
Green and Lazarsfeld.  We recall here these results in the  form which is more
convenient for our later use.  Next we deduce two statements about the holomorphic
Euler--Poincar\'e characteristic of a variety of maximal Albanese dimension fibered over
a curve.
\medskip

\begin{thm}\label{K2}
Let $X$, $Y$, $Z$ be projective varieties, $X$ smooth,   let  $f\colon X \to Y$, $g\colon Y\to Z$
be  surjective morphisms and let
$P_0\in
\Pic^0(X)$ be a torsion point. Then:
\begin{itemize}
\item[a)] $R^ig_*R^jf_*(\omega_X\ot P_0)$ is torsion free;
\item[b)] $R^ig_*R^jf_*(\omega_X\ot P_0)=0$ if $i>\dim Y-\dim Z$;
\item[c)] $R^k(g\circ f)_*(\ox\ot P_0)\cong \sum_iR^ig_*R^{k-i}f_*(\ox\ot P_0)$ ($\cong$ denotes
equality in the derived category).
\end{itemize}
\end{thm}
\begin{proof} If $P_0=0$ this is part of \cite[Thm. 3.4]{Ko2}. 
\newline If $P_0\ne 0$, then  let
$d$ be the order of $P_0$ in $\Pic^0(X)$,   let
$p\colon\tilde{X}\to X$ be the \'etale $\Z_d-$cover given by $P_0$ and let $\tilde{f}:=f\circ p$.
We have $p_*\omega_{\tilde{X}}=\oplus_{m=0}^{ d-1}\ \omega_X(mP_0)$,  hence
$R^jf_*(\omega_X\ot P_0)$ is a direct summand of $R^j\tilde{f}_*\omega_{\tilde{X}}$,
and  statements  a) and b) for $f$ follow by the corresponding statements for
$\tilde{f}$.

Consider now statement c). Since the morphism $p$ is finite, one has:
\begin{equation}\label{Rk1}
R^k(g\circ \tilde{f})_*(\omega_{\tilde{X}})=\oplus_{m=0}^{ d-1}R^k(g\circ
f)_*(\omega_X\ot mP_0)
\end{equation}
and the summands on the right hand side are the eigensheaves for the 
action of the Galois group $G\cong \Z_d$ of $p\colon \tilde{X}\to X$. 

On the other hand, by applying the case $P_0=0$ of c) one gets:
\begin{eqnarray}\label{Rk2}
R^k(g\circ \tilde{f})_*(\omega_{\tilde{X}})\cong
\sum_iR^ig_*R^{k-i}\tilde{f}_*(\omega_{\tilde{X}})=\\
= \sum_i\oplus_{m=0}^{d-1}\,R^ig_*R^{k-i}f_*(\ox\ot mP_0) .\nonumber
\end{eqnarray}
The result now follows by comparing the formulae \ref{Rk1} and \ref{Rk2}, provided that the
identification in  \ref{Rk2} is compatible with the $G-$action. This can indeed  be achieved by
choosing  the smooth ample divisor $H$ in the proof of \cite[Thm. 3.1]{Ko2} to be
$G-$invariant. 
\end{proof}

\medskip

 Let $A$ be an abelian variety and
$a\colon X\lra A$ a morphism. Given a sheaf $\F$ on $A$,  we define the cohomological support loci
$V^i(A,\F ):=\{P\in
\Pic^0(A)\,|\,h^i(\F\ot a^* P)>0\}$. When $a\colon X\lra \Alb(X)$ is the Albanese map, we
write simply
$V^i(\F)$.
\begin{thm}[Generic Vanishing Theorem]
\label{GVT}
\hfill\phantom{1}\par
\noindent Let $f\colon X\lra  Y$, $g\colon Y \lra  Z$ be  morphisms of smooth projective
varieties  and let $a\colon Z\lra  A$ be a morphism to an abelian
variety. Let 
$P_0\in \Pic ^0(X)_{tors}$ and set $\F :=R^kg_*R^jf_*(\omega _X \ot P_0)$. 
Then:
\begin{itemize}
\item[a)] $V^i(A,\F )$ has  codimension
$\ge i-(\dim Z-\dim a(Z))$;

\item[b)] Every irreducible component of $V^i(A,\F )$ is
a translate of a sub-torus of\, $\Pic^0(A)$ by a torsion point.\newline
The same also holds for the irreducible components of
$V^i_m(A, \FF):=\{ P\in \Pic ^0(A)\ |\ h^i(\FF \ot a^*P)\geq m\} $;

\item[c)] Let   $T$ be  an irreducible component of $V^i(A,\F )$,  let
$P \in T$ be a point such that $V^i(A,\F )$ is smooth at $P$, and let
$v\in H^1(A, \OO
_A)\cong T_{P }\Pic ^0(A)$. If $v$ is not tangent to $T$, then the
sequence
$$H^{i-1}(Y,\F \otimes a^* P) \stackrel {\cup v}  {\longrightarrow}
H^{i}(Y,\F \otimes a^*P)  \stackrel {\cup v}  {\longrightarrow }
H^{i+1}(Y, \F \otimes a^*P) $$
is exact. Moreover, if $P$ is a general point of $T$
and $v$ is
tangent to $T$ then
both maps  vanish;

\item[d)] If $Z$ has maximal Albanese dimension, then there are inclusions:
$$ V^0(A,\F )\supseteq V^1(A,\F)
\supseteq \dots \supseteq V^n(A, \F ).$$

\end{itemize}
\end{thm}
\begin{proof}Assume first that $X=Y=Z$ and that $f$ and $g$ are the identity map. 
If   $A$
is the Albanese variety and $P_0=0$, then this is the classical generic vanishing
theorem, which is due to Green--Lazarsfeld (\cite{GL1},
\cite{GL2}),  apart from the statement in b) that the components of $V^i(A,\F )$ are
translates by a torsion point, which is due to Simpson (\cite{Si}). The extension to the
case of a morphism
$a\colon X\lra A$ into any abelian variety is straightforward (cf. \cite[Rem. 1.6]{EL}).
To prove 
the statement in the case $P_0\ne 0$, 
one considers the 
\'etale $ \Z_d-$cover $p\colon \tilde{X}\lra X$ associated to $P_0$ as in the proof of Theorem
\ref{K2}. Then for $P\in\Pic^0(A)$ one has: $$h^i(\tilde{X},\omega_{\tilde{X}}\ot (a\circ
p)^*P)=\sum_{j=0}^{d-1} h^i(X,\omega_X\ot jP_0\ot a^*P).$$ Since $h^i(X,\omega_X\ot jP_0\ot
a^*P)$ is an upper semicontinuous function of $P\in \Pic^0(A)$, it is easy to show that if $T$ is
a component of $V^i_m(A,
\ox\ot P_0)$ then
$T$ is also a component of $V^i_{m'}(A, \omega_{\tilde{X}})$ for some $m'$.  Hence statements
a), b), c), d) follow from the corresponding statements for
the loci $V^i_m(A,\omega_{\tilde{X}})$. 

Consider now the general case.
By Theorem \ref{K2}, c), for $P\in \Pic^0(Z)$ one
has:
\begin{gather}
H^i(X,\ox\ot P_0\ot  (a\circ g\circ f)^*P )\!=\!\oplus_b H^b(Z, R^{i-b}(g\circ
f)_*(\ox\ot P_0)\ot a^*P)\!=\nonumber\\= \oplus_{b, c}H^b(Z,R^cg_*R^{i-b-c}f_*(\ox\ot
P_0)\ot a^*P).\nonumber
\end{gather}

In view  of this equality  and of the semicontinuity of the function
$h^b(Z, R^cg_*R^{i-b-c}f_*(\ox\ot P_0)\ot a^*P)$,
statements b), c) follow by the analogous  statements in the
case
$f=g=Id$  by the same argument as above. 
For the statements a), d), we refer to \cite{ClH} \S 10, \cite{Hac} Corollary 4.2
and \cite{ChH2} \S 2.

\end{proof}
\begin{cor}\label{Rzero}
 Let $X$ be a smooth variety of maximal Albanese dimension,
let
$f\colon X\lra Y$ be a surjective map and let $F$ be a general fiber of $f$. If $P_0\in
\Pic^0(X)_{tors}$ is general, then
$f_*(\ox\ot P_0)$ is torsion free of rank $\chi(\omega_F)$ and $R^if_*(\ox\ot P_0)=0$ for $i>0$.
\end{cor}
\begin{proof}
 Since $X$ is of maximal Albanese dimension, the restriction
 of the Albanese map $F\lra A(X)$ is generically finite onto its image.
 Since  $P_0$ is general, by Theorem \ref{GVT} we have
$h^i(\omega_F\ot P_0)=0$ for $i>0$ and $h^0(\omega_F\ot
P_0)=\chi(\omega_F\ot P_0)=\chi(\omega_F)$. It follows that for $i>0$,
the sheaf
$R^if_*(\ox\ot P_0)$ is supported on a proper closed subset of
$Y$ and that $f_*(\ox\ot P_0)$ has rank equal to $\chi(\omega_F)$. The claim now
follows since by Theorem
\ref{K2} the sheaves 
$R^if_*(\ox\ot P_0)$ are torsion free for every $i\ge 0$.
\end{proof}

\begin{thm} \label{chi}
Let $X$ be a smooth variety of maximal Albanese dimension  and let
$f\colon X\longrightarrow B$ be a  map  onto a smooth curve of
genus $\ge 2$. Then for general
$P_0\in \Pic ^0(X)_{tors}$, one has: 
$$\chi (\omega _X)= \chi (\omega _F)\chi (\omega _B)+
\deg (f_*(\omega _{X/B} \ot P_0))\geq \chi (\omega _F)\chi (\omega _B).$$
\end{thm}
\begin{proof}Let $P_0\in \Pic ^0(X)$ be a general torsion point
and $Q\in \Pic ^0(B)$ be a general point. Then Theorem \ref{GVT} implies:
$$\chi (\omega _X)=\chi (\omega _X\ot P_0\ot f^*Q)=h^0(\omega _X\ot P_0\ot
f^*Q)=h^0(f_*(\omega_X\ot P_0)\ot Q).$$ In turn, generic vanishing for the vector bundle
$f_*(\omega _{X}\ot P_0)$ on $B$ (cf. Theorem \ref{GVT}) implies:
$$h^0(f_*(\omega _X\ot P_0)\ot Q)=\chi (f_*(\omega _X\ot P_0)\ot Q)=
\chi (f_*(\omega _X\ot P_0)).$$ 
Let $F$ be  a general fiber of $f$. By Corollary \ref{Rzero}, $f_*(\ox\ot P_0)$ is a vector
bundle on $B$ of rank $\chi(\omega_F)$.
Now, since  $f_*(\omega _X\ot P_0)=\omega _B \ot f_*(\omega _{X/B}\ot P_0)$, 
 Riemann-Roch for vector bundles on a curve gives 
\begin{gather}
\chi (f_*(\omega _X\ot P_0))= \deg f_*(\ox\ot P_0)+(1-g(B))\chi(\omega_F)=
\nonumber\\= 
\deg (f_*(\omega _{X/B}\ot P_0))+(g(B)-1)\chi(\omega_F).\nonumber
\end{gather}
To finish the proof it is enough to show that $\deg(f_*(\omega _{X/B}\ot P_0))\geq 0$ for $P_0\in
\Pic^0(X)_{tors}$.  
If $P_0=0$, this follows from \cite{V} Satz V.
The analogous statement for any torsion point
$P_0$ follows from the case
$P_0=0$ by the same argument as in the proof of Theorem \ref{K2}.
\end{proof}

\noindent {\bf Remark:} It is not the case that for 
$f\colon X\lra Y$ a surjective morphism 
of varieties of maximal Albanese
dimension, one has $\chi (\omega _X)\geq \chi (\omega _Y)\chi (\omega _F)$.

Consider in fact for $i=1,2,3$ curves $C_i$ with elliptic involutions
$\sigma _i$. Then there is a morphism
$$X:=(C_1\times C_2 \times C_3)/\Z _2 \lra (C_2\times C_3)/\Z_2=:Y$$
where the groups $\Z _2$ act diagonally.
If $g(C_i)\geq 2$ for all $i$, it follows that $X$ is birational to
a smooth projective variety of general type
and $\chi (\omega _X)=0$ (cf. \cite{EL} or the computation in Example \ref{ex1}).
However $Y$ is a surface of general type and hence $\chi (\omega _Y)>0$.
\begin{prop}\label{chi0}
 Let $X$ be a smooth  variety of maximal Albanese 
dimension, $f\colon X\lra Y$ a surjective morphism and $F$ a general fiber of $f$.
If $\chi (\omega _F)=0$, then $\chi (\omega _X)=0$.
\end{prop}
\begin{proof} 
Let $P_0\in \Pic^0(X)_{tors}$ be a general point. By Corollary \ref{Rzero}, we have
$R^if_*(\ox\ot P_0)=0$ for every $i\ge 0$. It follows: $$\chi(\ox)=\chi(\ox\ot
P_0)=\sum_i(-1)^i\chi(R^if_*(\ox\ot P_0))=0.$$
\end{proof}

\section{The theorem}

This section is devoted to the proof of the following theorem:

\begin{thm}\label{mainT}
Let $X$ be a smooth projective variety of 
maximal Albanese dimension. If $\chi (\omega _X)=1$, then $q(X)\leq
2\dim X$ and if $q(X)=2\dim X$, then $X$ is birational to a product 
of curves of genus $2$.
\end{thm}

Our  proof uses induction on the dimension of $X$. The inductive step consists in 
showing that
$X$ admits a fibration $f\colon X\lra C$ onto a curve of genus $\ge 2$  and that the fibers of
$f$ satisfy the assumptions of Theorem \ref{mainT}. The existence of the
fibration $f$ is  in turn established  by  induction, exploiting the geometry of the
cohomological support loci (cf. Section \ref{pre}) of a suitable sheaf. We start by giving
an upper bound on the codimension of such a locus.

\begin{lem}\label{gut}
Let $f\colon X\lra Y$ be a morphism of smooth varieties of maximal 
Albanese dimension, let $P_0\in \Pic ^0(X)_{tors}$ be general, 
and  let $T$ be a component of $V^1(f_*(\omega _X \ot P_0))$. 
Assume that $\chi (\omega _X)=1$.

If  $T\not\subset V^i(f_*(\omega _X \ot P_0))$ then $\cod T\le i$.
\end{lem}
\begin{proof}
We write $\F :=f_*(\omega _X \ot P_0)$.
Since $P_0$ is general, by Corollary \ref{Rzero} we have 
$R^if_*(\ox\ot P_0)=0$ for $i>0$, and thus  $\chi(\F)=\chi(\ox\ot P_0)=1$. 

Let
$P\in T$ be a general  point,  so that $P$ is smooth for $V^1(\F)$ and  $P\notin V^i(\F)$. 
Assume
by contradiction that $T$ has codimension $\ge i+1$ and let
$W\subseteq H^1(\OO_X)$ be a linear subspace of dimension $i+1$ transversal to $T$ at $P$. 
Writing 
$\pp:=\pp(W)$,  the derivative  complexes of  Thm. \ref{GVT}, d) fit together to give an
exact complex of vector bundles (cf. \cite{EL}, proof of Theorem 3): 
\begin{equation}\label{complex}
0\!\to\!\OO_{\pp}(-d)\!\to\! H^0(\F\ot P)\ot \OO_{\pp}(-i)\!\to\!  \dots\! \to\!
H^{i-1}(\F\ot P)\ot\OO_{\pp}(-1)\!\to\! 0.
\end{equation}  Notice that we have $d\ge i$ by the exactness of the
complex.  Considering the hypercohomology spectral sequence
associated to the above complex, 
one sees that $H^{i}(\OO_{\pp}(-d))=0$, hence $d=i$.
Now twisting (\ref{complex}) by $\OO _{\pp}(-1)$ 
and considering the corresponding spectral  sequence,
one obtains an isomorphism $H^i(\OO_{\pp}(-i-1))\cong H^0(\F\ot P)\ot H^i(\OO_{\pp}(-i-1))$.
It follows that 
$h^0(\F\ot P)=1$ and that  the map
$H^{1}(\F\ot P)\ot\OO_{\pp}(-i+1)\lra H^2(\F\ot P)\ot\OO_{\pp}(-i+2)$ is injective.
Hence the following  complex is exact:
\begin{equation}
0\to H^1(\F\ot P)\ot \OO_{\pp}(-i-1)\to  \dots \to
H^{i-1}(\F\ot P)\ot\OO_{\pp}(-3)\to 0
\end{equation}
 By applying the previous argument again, we get $H^1(\F\ot P)=0$, a contradiction.  
\end{proof}
\begin{lem}\label{reldim}
Let $X$, $Y$ be smooth varieties, with $Y$ of maximal Albanese dimension, and
let $f\colon X\lra Y$ be a surjective morphism. Let $P_0\in \Pic^0(X)$ be a torsion point,
let $T$ be a component of $V^i(f_*(\omega_X\ot P_0))$ and let $g\colon Y\lra S:=\Pic^0(T)$
be the induced map. Then:
$$\dim g(Y)\le \dim Y-i.$$
\end{lem}
\begin{proof}
Since  $f_*(\omega_X\ot P_0)$ satisfies Theorem
\ref{GVT},   we have $T=Q_0+T_0$, where $Q_0\in
\Pic^0(Y)$ is a torsion point and $T_0\subset \Pic^0(Y)$ is a subtorous.
For $P\in T_0$ general, for $j\ge 0$ and  for $s>0$, by Theorem \ref{GVT} one has:
$$H^s(R^jg_*(f_*(\omega_X\ot P_0)\ot Q_0)\ot P)\!=\!H^s(R^jg_*(f_*(\omega_X\ot P_0\ot
f^*Q_0))\ot P)\!=\!0.$$ Hence for $P\in
T_0$ general, the Leray spectral sequence gives:
$$H^i(f_*(\omega_X\ot P_0)\ot P\ot Q_0)=H^0(R^ig_*(f_*(\ox\ot P_0\ot f^* Q_0))\ot P).$$ 
Now assume by contradiction that the relative dimension of $g$ is strictly smaller than $i$. Then
by Koll\'ar's relative vanishing (cf. Theorem \ref{K2}) we have $R^ig_*(f_*(\omega_X\ot
P_0\ot f^*Q_0))=0$, hence $H^i(f_*(\omega_X\ot P_0)\ot P\ot Q_0)=0$, contradicting the assumption
that
$T=T_0+Q_0$ is contained in $V^i(\F)$.

\end{proof}

\begin{prop} \label{P4}
Let $Y$ be a smooth variety  of maximal Albanese dimension such
that 
$q(Y)\ge  2\!\dim Y$. 
Assume that there exists a surjective map $f\colon X\lra Y$, with $X$ a smooth
variety of maximal Albanese dimension such that  $\chi(\omega_X)=1$.

Then there exists a surjective map $\psi\colon Y\lra C$, with $C$ a smooth curve of genus $\ge 2$.
\end{prop}
\begin{proof}
We prove the statement by induction on the dimension of $Y$, the case $\dim Y=1$ being obvious. 

Let $P_0\in \Pic^0(X)$ be a general torsion point and set $\F:=f_*(\omega_X\ot P_0)$. 
Since $P_0$ is general, we have  $R^if_*(\omega_X\ot P_0)=0$ for $i>0$ by Corollary \ref{Rzero}.
Therefore, we have $H^k(Y,
\F\ot P)=H^k(X, \omega_X\ot P_0\ot f^*P)$ for every $P\in \Pic^0(Y)$ and for every $k$. In
particular, 
$\chi(\F)=\chi(\omega_X\ot P_0)=1$.

Assume that $\dim V^1(\F)\geq 0$,  let $T\subseteq V^1(\F)$ be a component and let $i\ge 1$ be such
that $T\subseteq V^i(\F)$ and $T\not\subset V^{i+1}(\F)$. By Theorem \ref{GVT}, $T$ is a proper
subset of
$\Pic^0(Y)$. By Lemma
\ref{gut} we have
$\cod T\le i+1$. Let $g\colon Y\lra S:=\Pic^0(T)$ be the induced map and let $Z$ be a
nonsingular model of
$g(Y)$. Since $g(Y)$ generates $S$, we have 
$$q(Z)\ge \dim S\geq q(Y)-i-1.$$ By Lemma \ref{reldim},
we have $$\dim Z=\dim g(Y)\le \dim Y
-i.$$ 
Since $q(Y)\ge 2\!\dim Y$, it follows that 
$$q(Z)\ge 2\!\dim Y -i-1\ge 2\!\dim Z+i-1\geq 2\!\dim Z.$$ Hence by
the inductive hypothesis there exists a smooth curve $C$ of genus at
least $2$ and a surjective morphism
$h\colon Z\lra C$. We let $\psi\colon Y\lra C$ be the map obtained by composing the rational
map $Y\lra Z$ with $h$. 
In principle $\psi$ is just a rational map, but it is not difficult to show that $\psi$
is indeed regular. In fact, let $Y'\lra Y$ be a birational modification such that $\psi$ induces a
regular map $\psi'\colon Y'\lra C$ and let $\psi'_*\colon \Alb(Y)\lra J(C)$ be the induced
morphism. Then $\psi$ is the composition of the Albanese map $Y\lra \Alb(Y)=\Alb(Y')$ with
$\psi'_*$, hence it is a morphism.

To finish the proof we have to rule out the possibility that
 $V^1(\F)=\emptyset$. 
Assume that this is the case. 
Let $a\colon Y\lra A:=\Alb (Y)$ be the Albanese morphism and set $\F ':=a_* \F$. The map
$a$  is generically finite,  hence by Theorem
\ref{K2} we have $R^ia_* \F =0$ for all $i>0$. It follows that 
 $h^i(\F \ot a^* P)=h^i(\F '\ot P)$ for all $i\geq 0$ and $P\in \Pic ^0(Y)$.
Since $Y$ has maximal Albanese dimension, by Theorem \ref{GVT} the assumption
$V^1(\F)=\emptyset$ implies 
$h^i(\F '\ot P)=0$ for all
$i>0$ and any
$P\in
\Pic ^0(Y)$, and therefore $h^0(\F '\ot P)=\chi (\F'\ot P)=1$ for all $P\in \Pic ^0(Y)$.
Let $g=q(Y)$, $\hat {S}\colon D(A)\lra D(\Pic ^0(Y))$ be the Fourier--Mukai
transform on $A$ and $S\colon D(\Pic ^0(Y))\lra D(A)$ the Fourier--Mukai transform on
$\Pic^0(Y)$, so that by \cite{M} one has:
$$S\circ \hat {S}\cong (-1_A)^*[-g].$$
The sheaf $\hat {S} (\F ')=\hat {S} ^0(\F ')$ is a negative semidefinite 
line bundle $L$
on $\Pic ^0(Y)$. Since $L^\vee $ is semi-positive, there exist an abelian variety
$A'$, a surjective morphism $\Pic ^0(Y) \lra A'$ and a polarization $M$ on $A'$
such that $L$ is algebraically equivalent to the pull back of $M^\vee$ on $\Pic ^0(Y)$.
It follows that the support of $S^g(L)$ is
a translate of $\Pic ^0(A') \subset \Pic ^0( \Pic ^0 (A))=A$.  
On the other hand, we have that 
$S^g(L)=S^g(\hat {S}^0 (\F '))=(-1_A)^*(\F ')$
and by Theorem \ref{K2} the support of $\F '$ is equal to $a(Y)$. 
Since $a(Y)$ is a proper subvariety that 
generates $A$, we have obtained a contradiction and the proof is complete.
\end{proof}
We can  now give the proof of Theorem \ref{mainT}.
\begin{proof}[Proof of Theorem \ref{mainT}]We proceed by induction on $\dim X$.
By Proposition \ref{P4} (with $f=id$), there exists a surjective morphism
$\psi \colon X\lra C$ where $C$ is a smooth curve of genus at least $2$. Up to 
replacing
$\psi$ by the first term of its  Stein factorization, we may assume that
$\psi$ has connected fibers.
Let $F$ be a general fiber. Since $X$ is of general type and
maximal Albanese dimension, then $F$ is also
of general type and maximal Albanese dimension.
By Proposition \ref{chi0}, we have $\chi (\omega _F)>0$.
By Theorem \ref{chi}, one has
$$1=\chi (\omega _X)\geq \chi (\omega _F)\chi (\omega _C)>0$$
and so $ \chi (\omega _F)=\chi (\omega _C)=1$.
 The inductive assumption gives $q(F)\le 2\dim F
=2(\dim X-1)$ and therefore one has $q(X)\le q(F)+g(C)\le 2\dim X$.

Assume now that $q(X)=2\dim X=2n$. In this case one has $q(F)\ge q(X)-2\ge 2\dim F$, hence by
induction
$F$ is birational to a product $C_1\times ... \times C_{n-1}$ of curves of genus $2$. We need to
show that the fibration $\psi$ is  birationally a product.

The $3-$canonical map of $F$  is a birational morphism whose  
image is  isomorphic to $C_1\times ... \times C_{n-1}$.
Hence,  by replacing
$X$ with a desingularization of the relative $3-$canonical 
image, we may assume that the general fiber $F$
is isomorphic to a product $C_1\times ... \times C_{n-1}$. 
Notice that the decomposition $F=C_1\times ... \times C_{n-1}$ is unique,
up to permuting the $C_i$'s. Indeed, the curves of the form 
$$x_1\times ....\times x_{i-1}\times C_i\times x_{i+1}
\times ... \times x_{n-1}$$ are the only curves of genus 2 on 
$F$ with trivial normal bundle. Hence  every
automorphism of
$F$ either  preserves these systems of curves or permutes them 
(we have $C_i=C_{\sigma (i)}$ in this case, where $\sigma \in S_{n-1}$
is the corresponding permutation).
Using this remark it is easy to check also that the action of $\Aut(F)$ on $H^0(\Omega^1_F)$ is
faithful.

Denote by $P$ the cokernel of the injection $J(C)\lra \Pic^0(X)$.  Arguing as in the proof of the 
Lemme on page 345 of
\cite{Be}, one shows that there exists a nonempty open subset $C_0\subset C$ such that, denoting 
by $\psi_0\colon X_0\lra C_0$ the restriction of $\psi$, there is an isomorphism  
$\Pic^0(X_0/C_0)\overset{\sim}{\lra} C_0\times P$ over $C_0$. By the discussion above, for every $t\in
C_0$ there is a canonical (up to the order) splitting $\Pic^0(F_t)=J(C_{1,t})\times 
.. \times J(C_{n-1,t})$, where $F_t$  
denotes, as usual, the fiber over $t$. Since the abelian subvarieties of a fixed abelian variety are a
discrete set, it follows that the induced decomposition
$P=J_1\times ... \times J_{n-1}$ is independent of $F_t$. 

Take a finite covering $C'\lra C$ such that the induced fibration $X'\lra C'$ 
has a
section  (we may assume  $X'$ smooth, by possibly replacing it with a desingularization). 
Then we can consider the composition of the  relative Albanese map
$\Alb(X'/C')\lra C'\times P^{\vee}$  with the projection $C'\times P^{\vee}\lra C'\times 
J_1$. Let $Y$ be
a smooth projective model of the image of this map. By construction, $Y$ is a 
surface fibered in 
curves of genus
$2$ and has relative irregularity $\ge 2$. By the Lemme on page 345 of \cite{Be}, $Y$ is
birationally a product  and, in particular,  the curve $C_1$ 
does not vary in moduli as $t$
varies in $C$. The analogous statement holds of course for $C_2,...,C_{n-1}$, hence two  general
fibers
$F$ of $f\colon X\lra C$  are isomorphic.

Since the automorphism group of $F$ is finite, one can argue as in \cite{Se} and prove that $X$
is birational to a quotient $(F\times B)/G$, where $G$ is a finite group  that acts faithfully on $F$
and $B$ and diagonally on $F\times B$, and such that  $B/G=C$. One has $2n= q(X)=q(B/G)+q(F/G)$, namely
$q(F/G)=2(n-1)$. Since, as we have remarked above,  the representation of
$\Aut(F)$ on
$H^0(\Omega^1_F)$ is faithful, it follows that the group $G$ is trivial and that $X$ is birational
to
$C_1\times ... \times C_{n-1}\times C$. 

\end{proof}

\section{Examples}\label{Sexamples}

We describe here two families of varieties with $\chi
(\omega _X)=1$ and unbounded irregularity.

The first example   shows that in
Beauville's result the assumption that ``$X$ be a surface 
of general type''  cannot be weakened to
$\kappa(X)\ge 1$. The second example shows that in dimension $\ge 3$ the statement of Theorem
\ref{mainT} does not hold under the weaker assumption that $X$ be of general type, hence the
analog of  Beauville's result is not true in higher dimension.

\begin{ex}\label{ex1}
{\rm  We construct a family of properly elliptic surfaces 
$Y$ with $\chi(\omega _Y)=1$ and arbitrary irregularity
$g\ge 2$.

Let $B\/$ be a hyperelliptic curve of genus $g\ge 2$ and let $M$ be a line
bundle of degree $1$ on $B$ such
that the linear series $|2M|$ is the $g^1_2$ of $B$.
Set $\cF:=\pp(\OO_B\oplus \OO_B(2M))$, and 
let $p\colon \cF\lra B$ be the corresponding projection. Let
$\si_{\infty}$ be the section at ``infinity'' and set 
$\si_0:=\si_{\infty}+2p^*M$.
One has: 
$$\si_{\infty}^2=-2, \quad\si_0^2=2,\quad\si_0\si_{\infty}=0, \quad
K_{\cF}=-2\si_{\infty}+p^*(K_B-2M).$$ 
It is easy to check that
$|\si_0|$ is a free linear system of dimension $2$ on
$\cF$. 

Let $D_0\in |3\si_0|$ be a smooth divisor,  set $L:=2\si_{\infty}+3p^*M$,
$D=D_0+\si_{\infty}$  and let
$\pi\colon Y\lra \cF$ be the double cover associated to the relation $2L\equiv 
D$. The surface $Y$ is
smooth, since $D$ is smooth,  and $p$ induces an elliptic fibration $f\colon Y
\lra B$. Notice that
the surface
$Y$ is properly elliptic, since $g\ge 2$.

The standard formulae for double covers give:
\begin{gather*} \chi(Y)=2\chi(\cF )+
L(K_{\cF}+L)/2=1,\\
p_g(Y)=h^0(\cF, K_{\cF}+L)=h^0(B,K_B+M)=g.
\end{gather*}
It follows that $q(Y)=g$ and  $f\colon Y\lra B$ is the
Albanese map.

 Notice that if $C_1,\dots ,C_{n-2}$ are curves of genus 2, then the
variety $X:=Y\times C_1\times \dots \times C_{n-2}$ has dimension $n$ and it satisfies: 
$$\chi(\omega_X)=1,\quad  q(X)=g+2n-4, \quad \kappa(X)=n-1.$$}
\end{ex}

\begin{ex}\label{ex2}
{\rm We construct a family of $3-$folds $Y$  of general type with $\chi(\omega_Y)=1$ and
arbitrarily large irregularity.
Let $D$ be a smooth plane quartic  and let $S$ be
the second symmetric product of $D$. The surface $S$ is smooth minimal
of general type with
$q(S)=p_g(S)=3$.  Let $\si$ be the involution of $S$ that maps a pair $\{P, Q\}\in
S$ to the pair $\{P',Q'\}$ such that $P+Q+P'+Q'$ is the section of $D$ with a line in
$\pp^2$.  It is well known (cf. for instance \cite[Ex. 3.3]{cpt} ) that the fixed locus of $\si$
consists of 28 isolated fixed points and that the quotient surface
$\Si:=S/\si$ is a nodal surface of general type with $p_g(\Si)=3$ and $q(\Si)=0$.

Let $C$ be a smooth curve of genus $g\ge 2$,  let $C'\lra C$ be a smooth double cover
branched on $2m$ points and let $\tau$ be the corresponding involution of $C'$. 
Set $Z:=(S\times C')/(\si\times \tau)$. The variety $Z$ is of general type,  since it
dominates
$\Si\times C$,  and it has terminal singularities at    the
images of the
$56m$ fixed points of
$\si\times \tau$. A desingularization  $Y$ of $Z$ can be obtained as follows:   blow up
$X$ at the fixed points of $\si\times \tau$ to obtain a $3-$fold $\hat{X}$ such that
$\si\times \tau$ induces an involution of $\hat{X}$ whose fixed locus is a divisor and
then take  $Y$ to be the quotient of $\hat{X}$ by this involution. This description
shows that for every $i\ge 0$
one can identify $H^0(Y, \Omega_Y^i)$ with the subspace
of $H^0(S\times C', \Omega_{S\times C'}^i)$ invariant under the involution. 
We
write $H^0(C, \omega_{C'})=H^0(C, \omega_C)\oplus V$, where $V$ is the subspace of
antiinvariant forms,  which has dimension $g+m-1$. Then we have identifications: 
\begin{gather*}
H^0(\omega_Y)=H^0(\omega_C)\ot H^0(\omega_S), \quad 
H^0(\Omega^1_Y)=H^0(\omega_C),\\H^0(\Omega_Y^2)= H^0(\Omega^1_S)\ot V\oplus
H^0(\omega_S) 
\end{gather*}
It follows that $\chi(\omega_Y)=g-3m-1$, $q(Y)=g$ and, in particular the map $Y\lra C$
is the Albanese map. For  $g=3m+2$ we have the required examples.

As in the previous example, one can obtain smooth $n-$dimensional varieties of general
type with $\chi(\omega _Y)=1$ and with arbitrarily large irregularity by taking the product of the
above
$3-$folds with $n-3$ curves of genus 2.
 }
\end{ex}

\vspace{1truecm}

\begin{minipage}{13cm}
\parbox[t]{6.5cm}{Christopher D. Hacon\\
Department of Mathematics\\
University of Utah\\
155 South 1400 East, JWB 233\\
Salt Lake City, Utah 84112-0090, USA\\
hacon@math.ucr.edu
 } \hfill
\parbox[t]{5.5cm}{Rita Pardini\\
Dipartimento di Matematica\\
Universit\'a di Pisa\\
Via Buonarroti, 2\\
56127 Pisa, Italy\\
pardini@dm.unipi.it}
\end{minipage}
\end{document}